\documentclass{amsart}

\newtheorem{theorem}{Theorem}[section]

\newtheorem{proposition}[theorem]{Proposition}
\newtheorem{corollary}[theorem]{Corollary}
\newtheorem{problem}[theorem]{Problem}
\newtheorem{problems}[theorem]{Problems}
\newtheorem{conjecture}[theorem]{Conjecture}

\theoremstyle{definition}
\newtheorem{definition}[theorem]{Definition}

\newtheorem{example}[theorem]{Example}

\usepackage{amscd,amssymb}

\begin{document}

\title[Computing with matrix invariants]
{Computing with matrix invariants}

\author[Vesselin Drensky]
{Vesselin Drensky}
\address{Institute of Mathematics and Informatics,
Bulgarian Academy of Sciences, 1113 Sofia, Bulgaria}
\email{drensky@math.bas.bg}
\thanks{This is an improved version of the talk of the author
given at the Antalya Algebra Days VII on May 21, 2005.
Partially supported by Grant MM-1106/2001 of the Bulgarian Foundation for
Scientific Research.}
\subjclass[2000] {Primary: 16R30; Secondary:
05E05, 05E10, 13A50, 13P10, 16P90, 16S15}
\keywords{trace rings, invariant theory, matrix concominants, defining relations,
Hilbert series, symmetric functions, Schur functions}

\begin{abstract}
We present an introduction to the theory of the invariants under
the action of $GL_n({\mathbb C})$ by simultaneous conjugation of
$d$ matrices of size $n\times n$. Then we survey some results, old
or recent, obtained by a dozen of mathematicians, on minimal sets
of generators, the defining relations of the algebras of
invariants and on the multiplicities of the Hilbert series of
these algebras. The picture is completely understood only in the
case $n=2$. Besides, explicit minimal sets of generators are known
for $n=3$ and any $d$ and for $n=4$, $d=2$. The multiplicities of
the Hilbert series are obtained only for $n=3,4$ and $d=2$. For
$n>2$ most of the concrete results are obtained with essential use
of computers.
\end{abstract}
\maketitle

\section{Introduction to invariant theory}

All considerations in this paper are over an arbitrary field $K$
of characteristic 0. If not explicitly stated,
all vector spaces and algebras are over $K$, and the algebras
are unitary and commutative.

To get some idea about invariant theory,
we start with the following well known result from the undergraduate
algebra course.

Let $n\geq 2$ be an integer and let
\[
A=K[X]=K[x_1,\ldots,x_n]
\]
be the algebra of polynomials in $n$ variables. The algebra
of symmetric polynomials $A^{S_n}=K[X]^{S_n}$
consists of all polynomials $f(x_1,\ldots,x_n)\in A$
such that
\[
f(x_{\sigma(1)},\ldots,x_{\sigma(n)})=f(x_1,\ldots,x_n)
\]
for all permutations $\sigma$ in the symmetric group $S_n$.

\begin{theorem}
{\rm (i)} The algebra $A^{S_n}$ is generated by the elementary
symmetric functions
\[
e_1=x_1+\cdots+x_n=\sum x_i,
\]
\[
e_2=x_1x_2+x_1x_3+\cdots+x_{n-1}x_n=\sum_{i<j}x_ix_j,
\]
\[
\cdots
\]
\[
e_n=x_1\cdots x_n.
\]

{\rm (ii)} The presentation of the symmetric polynomials as polynomials of
$e_1,\ldots,e_n$ is unique, i.e.
\[
A^{S_n}\cong K[y_1,\ldots,y_n].
\]
\end{theorem}

As a measure {\it how many are the symmetric polynomials in $n$
variables} we may use the dimensions of the vector spaces
$A_k^{S_n}$ of the homogeneous symmetric polynomials of degree $k$.
For example, for $n=3$, the vector space $A_k^{S_3}$ has a basis
consisting of all products $e_1^{k_1}e_2^{k_2}e_3^{k_3}$, where
$k_1+2k_2+3k_3=k$. Instead, we consider the {\it generating function}
of the sequence $\text{\rm dim}(A_k^{S_n})$, namely
the {\it formal power series}
\[
H(A^{S_n},t)=\sum_{k\geq 0}\text{\rm dim}(A_k^{S_n})t^k,
\]
which is called the {\it Hilbert series} of $A^{S_n}$.
Since the set of products
\[
\{e_1^{k_1}e_2^{k_2}\cdots e_n^{k_n}\mid
k_1,\ldots,k_n\geq 0\}
\]
is a basis of $A^{S_n}$,
it is easy to see that
\[
H(A^{S_n},t)
=\sum_{k_j\geq 0}t^{k_1}t^{2k_2}\cdots t^{nk_n}
\]
\[
=(1+t+t^2+\cdots)(1+t^2+t^4+\cdots)\cdots(1+t^n+t^{2n}+\cdots)
=\prod_{i=1}^n\frac{1}{1-t^i}.
\]

Let us consider another example.

\begin{example}\label{cyclicly symmetric polynomials}
Let $n=3$, $A=K[x_1,x_2,x_3]$
and let $G$ be the subgroup of $S_3$ generated by
the cycle $\rho=(123)$. We denote by $A^G$ the algebra of
``{\it cyclicly symmetric polynomials}'',
\[
A^G=\{f(x_1,x_2,x_3)\in A \mid
f(x_2,x_3,x_1)=f(x_1,x_2,x_3)\}.
\]
One can show that $A^G$ is generated by the elementary symmetric functions
$e_1,e_2,e_3$ and one more polynomial
\[
f_4=x_1^2x_2+x_2^2x_3+x_3^2x_1,
\]
which is not symmetric. The sum
\[
a=(x_1^2x_2+x_2^2x_3+x_3^2x_1)+(x_1x_2^2+x_2x_3^2+x_3x_1^2)
\]
of $f_4$ and its ``other half''
$x_1x_2^2+x_2x_3^2+x_3x_1^2$,
is symmetric. The polynomial $b=f_4(a-f_4)$ is also symmetric.
Hence both $a,b$ belong to $K[e_1,e_2,e_3]$ and
$f_4$ satisfies the relation $f_4^2-af_4+b=0$ with coefficients
\[
a=\alpha(e_1,e_2,e_3),b=\beta(e_1,e_2,e_3)\in K[e_1,e_2,e_3].
\]
There exists a natural homomorphism
$\pi:K[y_1,y_2,y_3,y_4]\to A^G$ defined by
\[
y_1\to e_1,\quad y_2\to e_2,\quad y_3\to e_3,\quad
y_4\to f_4,
\]
and
\[
r(y_1,y_2,y_3,y_4)=y_4^2-\alpha(y_1,y_2,y_3)y_4+\beta(y_1,y_2,y_3)
\]
belongs to the kernel of $\pi$.
One can show that the ideal $\text{Ker}(\pi)$ is generated by
$r(y_1,y_2,y_3,y_4)$. Hence $A^G$ has a basis
\[
\{e_1^{k_1}e_2^{k_2}e_3^{k_3},\quad e_1^{k_1}e_2^{k_2}e_3^{k_3}f_4
\mid k_1,k_2,k_3\geq 0\}
\]
and, if $A^G_k$ is the vector space of the homogeneous cyclicly
symmetric polynomials of degree $k$, then the Hilbert series of $A^G$ is
\[
H(A^G,t)=\sum_{k\geq 0}\text{\rm dim}(A^G_k)t^k
\]
\[
=\sum_{k_i\geq 0}t^{k_1}(t^2)^{k_2}(t^3)^{k_3}
+\left(\sum_{k_i\geq 0}t^{k_1}(t^2)^{k_2}(t^3)^{k_3}\right)t^3
=\frac{1+t^3}{(1-t)(1-t^2)(1-t^3)}.
\]
\end{example}

Let $V_n$ be the vector space with basis $X=\{x_1,\ldots,x_n\}$ and let
the general linear group $GL_n=GL_n(K)=GL(V_n)$ act canonically on $V_n$.
If we identify $g\in GL_n$ with the invertible matrix
$g=(g_{ij})$, $g_{ij}\in K$, then the action of $g$ on $x_j$ is defined by
\[
g(x_j)=g_{1j}x_1+\cdots+g_{nj}x_n.
\]
This action is extended {\it diagonally} on $A=K[X]=K[x_1,\ldots,x_n]$:
\[
g(f(x_1,\ldots,x_n))=f(g(x_1),\ldots,g(x_n)),\quad g\in GL_n,\quad
f(x_1,\ldots,x_n)\in A.
\]

\begin{definition}
Let $G$ be any subgroup of $GL_n$.
(One says that the group $G$ is {\it linear}.)
The polynomial $f(x_1,\ldots,x_n)$ is called a $G$-{\it invariant} if
\[
g(f(x_1,\ldots,x_n))=f(x_1,\ldots,x_n)
\]
for all $g\in G$. We denote the algebra of $G$-invariants by $A^G$.
\end{definition}

For example, if $G$ is the subgroup of all permutational matrices
\[
g_{\sigma}=\sum_{j=1}^ne_{\sigma(j),j},\quad \sigma\in S_n,
\]
where the $e_{ij}$'s are the matrix units,
then $G\cong S_n$ and the algebra of $G$-invariants consists of all
symmetric polynomials in $n$ variables.

For $n=3$
and $G$ being the cyclic group of order 3 generated by the matrix
\[
g=\left(\begin{matrix}
0&0&1\\
1&0&0\\
0&1&0\\
\end{matrix}\right),
\]
the algebra of invariants $K[x_1,x_2,x_3]^G$ consists of all
cyclicly symmetric polynomials.

It is easy to see that the algebra of invariants is graded with respect
to the usual degree of polynomials. This means that if the polynomial
$f=f_0+f_1+\cdots+f_p$ is a $G$-invariant and $f_k$ is its
homogeneous component of degree $k$, then $f_k$ is also $G$-invariant.
If $K[X]^G_k$ is the homogeneous component of degree $k$ of $K[X]^G$, then
one may define the Hilbert series of $K[X]^G$ as the formal power series
\[
H(K[X]^G,t)=\sum_{k\geq 0}\text{\rm dim}(K[X]^G_k)t^k.
\]

We state of list of problems which are among
the main problems in invariant theory.
For more details, see some book on invariant theory, e.g.
Dieudonn\'e and Carrell \cite{DC} or Springer \cite{Sp};
for computational aspects of invariant theory see Sturmfels \cite{St}.

Recall that if $R$ is a finitely generated commutative algebra, then
$R$ is a homomorphic image of the polynomial algebra
with the same number of generators. Hence
$R\cong K[y_1,\ldots,y_m]/I$ for some $m$ and for some ideal
$I$. The equations $f(y_1,\ldots,y_m)=0$, $f\in I$,
are called {\it relations} of the algebra $R$
with respect to the given set of generators.
The elements of any generating set of the ideal $I$ are called
{\it defined relations} of $R$. A finitely generated algebra
with a finite system of defining relations is called
{\it finitely presented}.

\begin{problems}\label{main problems of invariant theory}
{\rm (i)} For a given linear group $G$, find a set of generators
of the algebra of $G$-invariants.

{\rm (ii)} Is it true that the algebra of invariants $K[X]^G$
is finitely generated for any subgroup $G$ of $GL_n$?

{\rm (iii)} If $K[X]^G$ is finitely generated, find its
defining relations.
Is $K[X]^G$ always finitely presented?

{\rm (iv)} When $K[X]^G$ is isomorphic to a polynomial algebra
(i.e. has no defining relations with respect to a suitable system of
generators)?

{\rm (v)} Calculate the Hilbert series of $K[X]^G$. Is it
a rational function? What kind of properties of $K[X]^G$ can be
recovered from its Hilbert series?
\end{problems}

We give short comments on the solutions of the above problems
in the case of finite groups.
Problem \ref{main problems of invariant theory} (i)
is solved by the following theorem of Emmy Noether.

\begin{theorem} {\rm (Endlichkeitssatz, \cite{No})}
The algebra of invariants $K[X]^G$ is finitely generated
for any finite subgroup $G$ of $GL_n$, by a set of invariants
of degree $\leq \vert G\vert$, the order of the group $G$.
\end{theorem}

It is well known that the linear operator of $K[X]$ (called the
{\it Reynolds operator})
\[
\rho: f(x_1,\ldots,x_n)\to\frac{1}{\vert G\vert}
\sum_{g\in G}g(f)
\]
is a projection of $K[X]$ onto $K[X]^G$, i.e. $\rho^2=\rho$
and $\rho(K[X])=K[X]^G$.
If we take $f=x_1^{k_1}\cdots x_n^{k_n}$ and
consider all monomials of degree $\leq \vert G\vert$, then
we obtain a system of generators of $K[X]^G$.

Problem \ref{main problems of invariant theory} (ii)
was the main motivation for the Hilbert 14-th Problem
from the famous list of 23 open problems presented by Hilbert
at the International Congress of Mathematicians in Paris, 1900,
\cite{H2}. Only in the late 1950's Nagata \cite{Na2} found
the first counterexamples to the 14-th Hilbert Problem.
We recommend the recent survey by Freudenburg \cite{Fr}
for the state-of-the-art of the problem and
for different approaches to construct other counterexamples.

\begin{theorem} {\rm (Hilbert Basissatz,
or Hilbert Basis Theorem, \cite{H1})}
Every ideal of $K[y_1,\ldots,y_m]$ is finitely generated.
\end{theorem}

As an immediate consequence of the Hilbert Basissatz
and the theorem of Emmy Noether
we obtain that for finite groups $G$ the algebra
$K[X]^G$ is finitely presented, i.e. the positive solution of
Problem \ref{main problems of invariant theory} (iii).
Since the degree of the generators of $K[X]^G$
is bounded, in principle, one can answer also the question which are
the defining relations of $K[X]^G$.

Problem \ref{main problems of invariant theory} (iv) is solved by
the following theorem of Chevalley-Shephard-Todd \cite{Ch, ST}.

\begin{theorem}
For a finite group $G$, the algebra of invariants $K[X]^G$
is isomorphic to a polynomial algebra
(in $n=\vert X\vert$ variables) if and only if $G\subset GL_n$
is generated by pseudo-reflections.
\end{theorem}

Recall that the matrix $g\in GL_n$ is called a {\it pseudo-reflection},
if it is similar to a matrix of the form
\[
\text{diag}(\xi,1,\ldots,1)=\left(\begin{matrix}
\xi&0&\cdots&0\\
0&1&\cdots&0\\
\vdots&\vdots&\ddots&\vdots\\
0&0&\cdots&1\\
\end{matrix}\right),
\]
it is a {\it reflection} if $\xi=-1$.
The symmetric group $S_n$ is generated by
the transpositions $\tau_j=(1j)$, $j=2,\ldots,n$, and
the images
$g_{\tau_j}=1-e_{11}-e_{jj}+e_{j1}+e_{1j}$
of $\tau_j$ under the embedding of $S_n$ into $GL_n$
are reflections. Hence $K[X]^{S_n}$ is isomorphic to a
polynomial algebra. On the other hand, the matrix $g$ from Example
\ref{cyclicly symmetric polynomials} (as well as its square
$g^2$) is similar to
\[
g=\left(\begin{matrix}
1&0&0\\
0&\varepsilon&0\\
0&0&\varepsilon^2\\
\end{matrix}\right),
\]
where $\varepsilon$ is a primitive third root of 1 in
some extension of $K$.
This implies that $G=\{1,g,g^2\}\subset GL_3$ cannot be generated by
pseudo-reflections and $K[x_1,x_2,x_3]^G$ is not isomorphic to
a polynomial algebra.

Finally, the following Molien formula \cite{Mo} answers
Problem \ref{main problems of invariant theory} (v).

\begin{theorem}
For any finite group $G\subset GL_n$
\[
H(K[X]^G,t)=\frac{1}{\vert G\vert}\sum_{g\in G}
\frac{1}{\text{\rm det}(1-gt)}.
\]
\end{theorem}

Here $\text{\rm det}(1-gt)$ is the determinant of the $n\times n$ matrix
$1-gt$. In particular, the Hilbert series $H(K[X]^G,t)$ is a rational
function. We refer to the book of Stanley \cite{St2} for the relationship
between the algebraic properties of graded algebras
and their Hilbert series.

The natural generalization of the case of finite linear groups
is the case of reductive groups. Recall that, by the Maschke theorem, the finite subgroups $G$
of $GL_n$ are completely reducible. This means that if, for a suitable basis
of $V_n$, the group $G$ consists of matrices of the block triangular form
\[
\left(\begin{matrix}
\ast&\ast\\
0&\ast\\
\end{matrix}\right),
\]
then, changing the basis we may assume that the matrices are of the
block diagonal form
\[
\left(\begin{matrix}
\ast&0\\
0&\ast\\
\end{matrix}\right).
\]
A similar property holds for the class of reductive subgroups of $GL_n$.
In particular, this holds for all classical groups ($GL_k$, $SL_k$,
$O_k$, $SO_k$, $Sp_k$, $U_k$).

\begin{theorem}
{\rm (i) (Hilbert-Serre, see e.g. \cite{DC})}
For any reductive subgroup $G$ of $GL_n$ the algebra of invariants
$K[X]^G$ is finitely generated.

{\rm (ii) (Molien-Weyl formula, see \cite{W1})}
If $G\subset GL_n({\mathbb C})$ is compact, then one can define
Haar measure on $G$, replace in the Molien formula the sum with
an integral and obtain the formula for the Hilbert series
of the algebra of invariants ${\mathbb C}[X]^G$.
\end{theorem}

\section{Matrix invariants}

In the sequel, we do not use the symbol $X$ for the set of variables
$\{x_1,\ldots,x_n\}$ and we denote by $X,X_i,Y,Y_i$, etc. $n\times n$
matrices.
We start with a special case of the main object of this paper. Let
\[
X=(x_{ij})=\left(\begin{matrix}
x_{11}&\cdots&x_{1n}\\
\vdots&\ddots&\vdots\\
x_{n1}&\cdots&x_{nn}\\
\end{matrix}\right)
\]
be an $n\times n$ matrix with
$n^2$ algebraically independent entries $x_{ij}$.
Such a matrix is called a {\it generic} $n\times n$ matrix
because every $n\times n$ matrix $A=(a_{ij})$
with entries $\alpha_{ij}$ from a commutative $K$-algebra
$S$ can be obtained from $X$ by specialization
of the variables $x_{ij}$ (i.e. a homomorphism
$K[x_{ij}]=K[x_{ij}\mid i,j=1,\ldots,n]\to S$ defined by
$x_{ij}\to\alpha_{ij}$, $i,j=1,\ldots,n$).
The group $GL_n$ acts on $X$ by conjugation:
\[
gXg^{-1}=Y=(y_{ij}),\quad g\in GL_n,
\]
where the $y_{ij}$'s are linear combinations of the $x_{ij}$'s.
We define an action of $GL_n$ on the algebra
$K[x_{ij}]$ by $g:x_{ij}\to y_{ij}$, $i,j=1,\ldots,n$,
and are interested in the algebra of invariants $K[x_{ij}]^{GL_n}$.
The obvious invariants come from the Cayley-Hamilton theorem.
If
\[
f_X(t)=\text{\rm det}(1-tX)=(-1)^n(t^n+a_1t^{n-1}+\cdots+a_{n-1}t+a_n)
\]
is the characteristic polynomial of $X$, then the coefficients
$a_1,\ldots,a_n$ depend on the entries of $X$. In particular,
$a_1=-\text{\rm tr}(X)$ and $a_n=(-1)^n\text{\rm det}(X)$.
The Cayley-Hamilton theorem states that $f_X(X)=0$.
Since $gf_X(X)g^{-1}=f_X(gXg^{-1})$ for any $g\in GL_n$, we obtain that
the coefficients $a_1,\ldots,a_n$ are $GL_n$-invariants.
The following theorem is a partial case of Theorem \ref{first fundamental theorem} below,
but probably was known much before it. It seems impossible to trace its origin.

\begin{theorem}\label{invariants of single matrix}
With respect to the above action of $GL_n$, the algebra of invariants
$K[x_{ij}]^{GL_n}$ is generated by the coefficients of the characteristic
polynomial of $X$.
\end{theorem}

Up to a sign, the coefficients of the characteristic polynomial
are equal to the elementary symmetric polynomials evaluated on the
eigenvalues of the matrix. Let
\[
p_k=x_1^k+\cdots+x_n^k
\]
be the $k$-th power symmetric function. Since the characteristic of $K$ is equal to 0,
the Newton formulas
\[
p_k-e_1p_{k-1}+\cdots+(-1)^{k-1}e_{k-1}p_1
+(-1)^kke_k=0,\quad k=1,\ldots,n,
\]
allow to express the elementary symmetric polynomials
$e_1,\ldots,e_n$ in terms of the power symmetric polynomials. Hence $K[x_1,\ldots,x_n]^{S_n}$ is generated by
$p_1,\ldots,p_n$. If $\xi_1,\ldots,\xi_n$ are the eigenvalues
of the matrix $X$ in the algebraic closure of the field of rational
functions $K(x_{ij})$, then $X$ is similar to the upper triangular matrix
\[
\left(\begin{matrix}
\xi_1&\ast&\cdots&\ast\\
0&\xi_2&\cdots&\ast\\
\vdots&\vdots&\ddots&\vdots\\
0&0&\cdots&\xi_n\\
\end{matrix}\right).
\]
(It is known that $X$ is similar even to the
diagonal matrix $\text{diag}(\xi_1,\ldots,\xi_n)$
with the same diagonal entries as the above matrix.)
The coefficients of the characteristic polynomial $f_X(X)$ are
\[
a_1=-e_1(\xi_1,\ldots,\xi_n),
a_2=e_2(\xi_1,\ldots,\xi_n),\ldots,
a_n=(-1)^ne_n(\xi_1,\ldots,\xi_n).
\]
Similarly,
\[
p_1(\xi_1,\ldots,\xi_n)=\text{\rm tr}(X),
p_2(\xi_1,\ldots,\xi_n)=\text{\rm tr}(X^2), \ldots,
p_n(\xi_1,\ldots,\xi_n)=\text{\rm tr}(X^n).
\]
In this way, as an immediately consequence of Theorem
\ref{invariants of single matrix}, we obtain:

\begin{corollary}
For the above action of $GL_n$ on $K[x_{ij}]$, the algebra
$K[x_{ij}]^{GL_n}$ is generated by the traces
$\text{\rm tr}(X),\text{\rm tr}(X^2),\ldots,\text{\rm tr}(X^n)$
of the powers of $X$.
\end{corollary}

Now we consider the more general situation, when the group $GL_n$
acts on several generic matrices. We fix an integer $d\geq 2$.
Let
\[
X_i=\left(x_{pq}^{(i)}\right),\quad p,q=1,\ldots,n,\quad i=1,\ldots,d,
\]
be $d$ generic $n\times n$ matrices and let
\[
\Omega=\Omega_{nd}=K[x_{pq}^{(i)}\mid p,q=1,\ldots,n,\quad i=1,\ldots,d]
\]
be the polynomial algebra in $dn^2$ variables. The action of $GL_n$ on $\Omega$
is as in the case of one generic matrix. If
\[
gX_ig^{-1}=Y_i=(y_{pq}^{(i)}),\quad g\in GL_n,
\]
then we define
\[
g:x_{pq}^{(i)}\to y_{pq}^{9i)},\quad p,q=1,\ldots,n,\quad i=1,\ldots,d.
\]
This action of $GL_n$ on $\Omega$ is called the action of $GL_n$
by {\it simultaneous conjugation} on $d$ generic $n\times n$ matrices.
Here are the main problems concerning the {\it algebra of matrix invariants}
$\Omega_{nd}^{GL_n}$.

\begin{problem}
Describe the algebraic properties of $\Omega_{nd}^{GL_n}$.
In particular, find concrete sets of generators and defining relations,
calculate the Hilbert series, etc.
\end{problem}

We refer e.g. to the books by Procesi \cite{P1}, Formanek \cite{F1},
and Drensky and Formanek \cite{DF} as a further reading on invariant theory
of matrices.

Since the action of $GL_n$ on $\Omega_{nd}$ is reductive, Classical
Invariant Theory gives:
\begin{theorem}
The algebra $\Omega_{nd}^{GL_n}$ is finitely generated.
\end{theorem}

Traditionally, a result giving the explicit generators of
the algebra of invariants of a linear group $G$ is called a {\it first
fundamental theorem of the invariant theory of} $G$ and a result
describing the relations between the generators is a {\it second
fundamental theorem}.
As in the case of a single generic matrix, the coefficients of
the Cayley-Hamilton polynomial of any matrix expressed
in terms of $X_1,\ldots,X_d$ is $GL_n$-invariant.
In particular, all traces $\text{\rm tr}(X_{i_1}\cdots X_{i_k})$
are $GL_n$-invariant. A more sophisticated result than the previous theorem
states:

\begin{theorem} {\rm (The First Fundamental Theorem of Matrix
Invariants)}\label{first fundamental theorem}
The algebra $\Omega_{nd}^{GL_n}$ is generated by all traces
$\text{\rm tr}(X_{i_1}\cdots X_{i_k})$, where $i_1,\ldots,i_k=1,\ldots,d$.
\end{theorem}

It is difficult to judge about the priority in the proof of this theorem.
(See the history of generic matrices in the paper by Formanek \cite{F2}.)
The proof follows from general facts on invariant theory of $GL_n$ and can
be found in the book by Gurevich \cite{Gu} and in the papers by Sibirskii \cite{Si}
and Procesi \cite{P2}. Although it was stated as well known by
Kirillov \cite{Ki},
it seems that the understanding of the importance of the theorem
is a result of its rediscovery by Procesi \cite{P2}.

Combining both theorems, we are looking for a concrete finite set of
generators of $\Omega_{nd}^{GL_n}$. We recall one of the most important
theorems in the theory of algebras with polynomial identities.

\begin{theorem} {\rm (Nagata-Higman Theorem)}
If the (nonunitary) algebra $R$ is nil of bounded index $\leq n$, i.e.
$r^n=0$ for all $r\in R$, then $R$ is nilpotent, i.e. there exists an
$N=N(n)$
such that $r_1\cdots r_N=0$ for all $r_1,\ldots,r_N\in R$.
\end{theorem}

The Nagata-Higman theorem was established in 1953 by Nagata \cite{Na1}
for algebras over a field of characteristic 0 and then in 1956
generalized by Higman \cite{Hi} when $\text{char K}=p>n$.
Much later it was discovered
that this theorem was first established in 1943 by Dubnov
and Ivanov \cite{DI} but their paper was overlooked
by the mathematical community.

The class of nilpotency $N(n)$ in the Nagata-Higman theorem is related
in the following nice way with invariant theory of matrices.

\begin{theorem}
{\rm (Formanek \cite{F11}, Procesi \cite{P2, P4}, Razmyslov \cite{R1})}
Let $N(n)$ be the class of nilpotency in the Nagata-Higman theorem.
Then the algebra of invariants $\Omega_{nd}^{GL_n}$ is generated
by the traces $\text{\rm tr}(X_{i_1}\cdots X_{i_m})$ of degree $\leq N(n)$.
For $d$ sufficiently large this bound is sharp.
\end{theorem}

It is important to know the
exact value of the class of nilpotency $N(n)$
in the Nagata-Higman theorem.
The upper bound given in the proof of Higman \cite{Hi} is
$N(n)\leq 2^n-1$.
The best known upper bound is due to Razmyslov \cite{R1}.
Applying {\it trace polynomial identities} of matrices,
he obtained the bound
$N(n)\leq n^2$. The proof of the theorem of Razmyslov is given
also in his book \cite{R2} or in the book by Formanek \cite{F1}.
For a lower bound, Kuzmin \cite{Ku} showed that
$N(n)\geq \frac{1}{2}n(n+1)$. A proof of the result of Kuzmin
may be found also in the books by Drensky and Formanek \cite{DF} or by Kanel-Belov and Rowen
\cite{KBR}.
Hence
\[
\frac{n(n+1)}{2}\leq N(n) \leq n^2.
\]

\begin{problem} Find the exact value $N(n)$ of the
class of nilpotency of nil algebras of index $n$
over a field of characteristic $0$.
\end{problem}

\begin{conjecture}
{\rm (Kuzmin \cite{Ku})} The
exact value $N(n)$ of the
class of nilpotency of nil algebras of index $n$
over a field of characteristic $0$ is
\[
N(n)=\frac{n(n+1)}{2}.
\]
\end{conjecture}

The only values of $N(n)$ are known for $n \leq 4$:
Dubnov \cite{Du} obtained in 1935
\[
N(1) = 1,\quad N(2) = 3,\quad N(3) = 6.
\]
In 1993, Vaughan-Lee \cite{VL} proved that
\[
N(4)=10.
\]
In this way the conjecture of Kuzmin is confirmed for $n\leq 4$.
Recently, Shestakov and Zhukavets \cite{SZ} have proved that the class
of nilpotency of the two-generated algebras satisfying the identity
$x^5=0$ is equal to 15, which agrees with the conjecture of Kuzmin for $n=5$.
They have obtained the same result also in the more general setup
of 2-generated superalgebras. Their proof is based on computer calculations
with the GAP package.

Since $N(2)=3$, the algebra $\Omega_{2d}^{GL_2}$ is generated by
products of traces of degree $\leq 3$. The following result was established
by Sibirskii \cite{Si}.

\begin{theorem}
The elements
\[
\text{\rm tr}(X_i),\text{\rm tr}(X_i^2),\quad i=1,\ldots,d,\quad
\text{\rm tr}(X_iX_j),\quad 1\leq i<j\leq d,
\]
\[
\text{\rm tr}(X_iX_jX_k),\quad 1\leq i<j<k\leq d,
\]
constitute a minimal set of generators of the
algebra of $2\times 2$ matrix invariants $\Omega_{2d}^{GL_2}$.
\end{theorem}

Now we give an idea about the Razmyslov-Procesi theory which is related with the second
fundamental theorem of the matrix invariants, see
Razmyslov \cite{R1} and Procesi \cite{P2}, as well as
the book by Razmyslov \cite{R2} for other applications
of his method. For simplicity we consider the case $n=2$ only.
The Cayley-Hamilton theorem for $2\times 2$ matrices implies that
\[
X^2 - \text{\rm tr}(X)X +\text{\rm det}(X)=0.
\]
The Newton formulas give that
\[
\text{\rm det}(X)=\frac{1}{2}(\text{\rm tr}^2(X)
-\text{\rm tr}(X^2)).
\]
This can be seen also directly. If $\xi_1,\xi_2$ are
the eigenvalues of $X$, then
\[
\text{\rm tr}(X)=\xi_1+\xi_2,\quad
\text{\rm tr}(X^2)=\xi_1^2+\xi_2^2,
\]
\[
\text{\rm det}(X)=\xi_1\xi_2=\frac{1}{2}((\xi_1+\xi_2)^2-(\xi_1^2+\xi_2^2))
=\frac{1}{2}(\text{\rm tr}^2(X)
-\text{\rm tr}(X^2)).
\]
In this way we obtain the {\it mixed trace identity}
\[
c(X) = X^2 - \text{\rm tr}(X)X + \frac{1}{2}(\text{\rm tr}^2(X)
-\text{\rm tr}(X^2)) = 0.
\]
Now we consider the identity
$c(X_1+X_2)-c(X_1)-c(X_2)=0$, i.e. we
{\it linearize} the identity $c(X) = 0$. In this way we obtain
the {\it mixed Cayley-Hamilton identity}
\[
\Psi_2(X_1,X_2) = X_1X_2 + X_2X_1 - \text{\rm tr}(X_1)X_2 - \text{\rm tr}(X_2)X_1
+ \text{\rm tr}(X_1)\text{\rm tr}(X_2) - \text{\rm tr}(X_1X_2) = 0.
\]
Since the trace is a nondegenerate bilinear form on $M_2(K)$, the
vanishing of the polynomial $\Psi_2(X_1,X_2)$ on $M_2(K)$ is equivalent
to the vanishing of the {\it pure Cayley-Hamilton identity}
\[
\Phi_2(X_1,X_2,X_3)=\text{\rm tr}(\Psi_2(X_1,X_2)X_3)=0
\]
on all $2 \times 2$ matrices. Direct calculations show that
\[
0 = \Phi_2(X_1,X_2,X_3)=\text{\rm tr}(\Psi_2(X_1,X_2)X_3)
= \text{\rm tr}(X_1X_2X_3) +
\text{\rm tr}(X_2X_1X_3)
\]
\[
- \text{\rm tr}(X_1)\text{\rm tr}(X_2X_3)
- \text{\rm tr}(X_2)\text{\rm tr}(X_1X_3) +
\text{\rm tr}(X_1)\text{\rm tr}(X_2)\text{\rm tr}(X_3) -
\text{\rm tr}(X_1X_2)\text{\rm tr}(X_3).
\]
If we delete the symbols of traces and the $X$'s in the above expression,
we shall obtain the following linear combination of permutations
\[
(123)+(213)-(1)(23)-(2)(13)+(1)(2)(3)-(12)(3)=\sum_{\sigma\in S_3}\text{sign}(\sigma).
\]
This suggests the following construction. We write the permutations in the symmetric group
$S_m$ as products of disjoint cycles, including the cycles of length 1,
\[
\sigma=(i_1 \ldots i_p)(j_1 \ldots j_q) \cdots (k_1\ldots k_r).
\]
We define the {\it associated trace function}
\[
\text{\rm tr}_{\sigma}(x_1, \ldots ,x_m)
= \text{\rm tr}(x_{i_1} \cdots x_{i_p})
\text{\rm tr}(x_{j_1} \cdots x_{j_q})\cdots \text{\rm tr}(x_{k_1} \cdots x_{k_r}).
\]
For every element
\[
\sum_{\sigma\in S_m}\alpha_{\sigma}\sigma\in KS_m,
\quad \alpha_{\sigma} \in K, \sigma\in S_m,
\]
where $KS_m$ is the group algebra of $S_m$, we define the trace polynomial
\[
f(x_1,\ldots,x_m) = \sum_{\sigma\in S_m}\alpha_{\sigma}
\text{\rm tr}_{\sigma}(x_1, \ldots ,x_m).
\]
We also assume that for $m \leq k$ the symmetric group $S_m$ acts on
$1,\ldots,m$ and leaves invariant $m + 1,\ldots,k$, i.e. $S_m$ is
canonically embedded into $S_k$.

\begin{theorem} {\rm (The Second Fundamental Theorem of Matrix
Invariants, Raz\-mys\-lov \cite{R1}, Procesi \cite{P2})} Let
\[
f(x_1,\ldots,x_m) = \sum_{\sigma\in S_m}\alpha_{\sigma}
\text{\rm tr}_{\sigma}(x_1, \ldots ,x_m),\quad \alpha_{\sigma} \in K,
\]
be a multilinear trace polynomial of degree $m$. Then $f = 0$ is
a trace identity for the $n \times n$ matrix algebra, i. e.
$f(a_1,\ldots,a_m) = 0$ for all $a_1,\ldots,a_m \in M_n(K)$, if and
only if
\[
\sum_{\sigma\in S_m}\alpha_{\sigma}\sigma
\]
belongs to the
two-sided ideal $J(n,m)$ of the group algebra $KS_m$ generated by the
element
\[
\sum_{\sigma\in S_{n+1}}\text{\rm sign}(\sigma)\sigma.
\]
\end{theorem}
As in the case of $2\times 2$ matrices,
the {\it fundamental trace identity}
\[
\sum_{\sigma\in S_{n+1}}(\text{\rm sign } \sigma)
\text{\rm tr}_{\sigma}(x_1, \ldots, x_{n+1})=0
\]
is actually the linearization of the Cayley-Hamilton polynomial.

There are several important objects related with invariant theory of matrices.
As above, $n,d\geq 2$ are fixed integers and $X_1,\ldots,X_d$
are $d$ generic $n\times n$ matrices:

The algebra $R_{nd}$ generated by $X_1,\ldots,X_d$;

The {\it pure} (or {\it commutative}) {\it trace algebra} $C_{nd}=\Omega_{nd}^{GL_n}$
generated by the traces of all products $\text{tr}(X_{i_1}\cdots X_{i_k})$;

The {\it mixed} (or {\it noncommutative}) {\it trace algebra} $T_{nd}$ generated by $R_{nd}$
and $C_{nd}$ regarding the elements of $C_{nd}$ as scalar matrices;

The field of fractions $Q(C_{nd})$ of the algebra $C_{nd}$;

The algebra $Q(C_{nd})R_{nd}$.

All these algebras have no zero divisors and play important roles in mathematics.
See the books \cite{DF, F1, J, P1} for different aspects
of the theory of algebras of matrix invariants and their applications to combinatorial and structure theory
of PI-algebras, central division algebras, etc.

The algebra $R_{nd}$ is a well known object in the theory of
{\it PI-algebras},
or {\it algebras with polynomial identities}.
Let $K\langle x_1,\ldots,x_d\rangle$
be the free associative algebra
(i.e. the algebra of polynomials in noncommuting variables).
A polynomial $f(x_1,\ldots,x_d)\in K\langle x_1,\ldots,x_d\rangle$ is a {\it polynomial
identity} in $d$ variables for $M_n(K)$ if $f(a_1,\ldots,a_d)=0$ for all $a_1,\ldots,a_d\in M_n(K)$.
The set $I(M_n(K))$ of all polynomial identities in $d$ variables is a two-sided ideal of
$K\langle x_1,\ldots,x_d\rangle$ and $R_{nd}$ is isomorphic to the factor algebra
$K\langle x_1,\ldots,x_d\rangle/I(M_n(K))$.

Clearly, the algebra $C_{nd}$ is the algebra of matrix
invariants.

The algebra $T_{nd}$ is also the algebra of invariant polynomial functions under a suitable
action of $GL_n$. It is called the algebra of {\it matrix concominants}.
{\it It is a finitely generated $C_{nd}$-module and, as a $C_{nd}$-module,
has a generating set consisting of products $X_{j_1}\cdots X_{j_k}$,
where $k < N(n)$, the class of nilpotency in the Nagata-Higman theorem}.

The field of fractions
$Q(C_{nd})$ appears naturally in field theory. One of the main problems related with $Q(C_{nd})$
is {\it whether it is a purely transcendent extension of $K$}.

Finally, $Q(C_{nd})R_{nd}$ is a central
division algebra of dimension $n^2$ over its centre $Q(C_{nd})$ and serves as a source of counterexamples
to the theory of central division algebras.

General invariant theory gives that $C_{nd}$ and $T_{nd}$ have nice algebraic properties.

\begin{theorem} {\rm (Van den Bergh, \cite{VB1})}
The algebra $C_{nd}$ is a Cohen-Macaulay and even Gorenstein unique factorization domain.
The algebra $T_{nd}$
is a Cohen-Macaulay module over $C_{nd}$.
\end{theorem}

Recall that the Noether normalization theorem gives that $C_{nd}$
contains a homogeneous set of algebraically independent elements $\{a_1,\ldots,a_k\}$,
where $k=(d-1)n^2+1$ is the transcendence degree of the quotient field $Q(C_{nd})$
($k$ is also equal to the Krull dimension of $C_{nd}$),
such that $C_{nd}$ is integral over the polynomial algebra $K[a_1,\ldots,a_k]$.
Such a set $\{a_1,\ldots,a_k\}$ is called a
{\it homogeneous system of parameters} for
$C_{nd}$. By a result of Stanley \cite{St1}, {\it a graded $C_{nd}$-module is
Cohen-Macaulay if and only if it is a free module with respect to some
homogeneous system of parameters $\{a_1,\ldots,a_k\}$ of $C_{nd}$}.

The algebras $R_{nd},C_{nd}$, and $T_{nd}$ are multigraded and their homogeneous components
of degree $(k_1,\ldots,k_d)$ consist of all polynomials which are homogeneous of degree
$k_i$ in the generic matrix $X_i$. Hence we may consider their Hilbert series in $d$ variables.
For example,
\[
H(C_{nd},t_1,\ldots,t_d)
=\sum\text{\rm dim}\left(C_{nd}^{(k_1,\ldots,k_d)}\right)t_1^{k_1}\cdots t_d^{k_d},
\]
where $C_{nd}^{(k_1,\ldots,k_d)}$ is the homogeneous component of degree $(k_1,\ldots,k_d)$.
Le Bruyn \cite{LB1} for the case $n = 2$, and Formanek \cite{F10} and Teranishi \cite{T1, T3}
in the general case proved:

\begin{theorem}
Let $d \geq 2$ for $n\geq 3$ and $d > 2$ for $n = 2$.
Then $H(C_{nd},t_1,\ldots,t_d)$ and $H(T_{nd},t_1,\ldots,t_d)$
satisfy the functional equation
\[
H(C_{nd},t_1^{-1},\ldots,t_d^{-1})
= (-1)^k(t_1 \cdots t_d)^{n^2}H(t_1,\ldots,t_d),
\]
where $k = (d - 1)n^2 + 1$
is the Krull dimension of $C_{nd}$, and similarly for the Hilbert series of $T_{nd}$.
\end{theorem}

The proofs of this theorem given by Formanek and
Teranishi are quite different and use, respectively, representation theory
of general linear groups and the Molien-Weyl integral formula.
Later, Van den Bergh paid attention that the proof can be considerably simplified
using results of Stanley on Hilbert series of Cohen-Macaulay algebras.

We need some background on symmetric polynomials and
representation theory of the general linear group,
see e.g. the books by Weyl \cite{W2} and Macdonald \cite{Mc}.
As in the case of polynomial algebras,
the general linear group $GL_d$ acts canonically on the vector space
with basis $\{X_1,\ldots,X_d\}$. If
$g=(g_{ij})$, $g_{ij}\in K$, then the action of $g$ on $X_j$ is defined by
\[
g(X_j)=g_{1j}X_1+\cdots+g_{dj}X_d.
\]
This action is extended diagonally on $R_{nd},C_{nd},T_{nd}$.
If $f(X_1,\ldots,X_d)$ is any polynomial expression
depending on $X_1,\ldots,X_d$ (maybe including also traces), then
\[
g(f(X_1,\ldots,X_d))=f(g(X_1),\ldots,g(X_d)),\quad g\in GL_d.
\]
Representation theory of $GL_d$ says that every submodule of
the $GL_d$-modules $R_{nd},C_{nd},T_{nd}$ is a direct sum of
irreducible (or simple) submodules.
The irreducible $GL_d$-submodules which appear in
the decomposition are {\it polynomial modules} and are described
in terms of {\it partitions} of integers. If
\[
\lambda=(\lambda_1,\ldots,\lambda_d),\quad
\lambda_1\geq\cdots\geq\lambda_d\geq 0,
\]
is a partition of $k$ (notation $\lambda\vdash k$)
in not more than $d$ parts, then we denote the related $GL_d$-module by
$W(\lambda)=W(\lambda_1,\ldots,\lambda_d)$. To be explicit, below
we consider the case of $C_{nd}$ only. If
\[
C_{nd}=\bigoplus m(\lambda)W(\lambda),
\quad m(\lambda)\geq 0,
\]
i.e. there are $m(\lambda)$ direct summands isomorphic to $W(\lambda)$,
then we say that $W(\lambda)$ appears with {\it multiplicity}
$m(\lambda)$. The multiplicities $m(\lambda)$
for $R_{nd}, C_{nd}$, and $T_{nd}$ play important role in the
quantitative study of polynomial identities of matrices. (See
the survey by Regev \cite{Re} and the book of the author \cite{D}
for applications of representation theory of $S_n$ and $GL_d$
to the theory of PI-algebras.)
The Hilbert series of $C_{nd}$ has the form
\[
H(C_{nd},t_1,\ldots,t_d)
=\sum m(\lambda)S_{\lambda}(t_1,\ldots,t_d),
\]
where $S_{\lambda}(t_1,\ldots,t_d)$ is the {\it Schur function}
associated to $\lambda$.
Schur functions are important combinatorial objects and
appear in many places in mathematics. For example, {\it they form a basis
of the vector space of all symmetric polynomials in $d$ variables}.
One of the possible ways to define Schur functions is via
Vandermonde-like determinants. For a partition
$\mu=(\mu_1,\ldots,\mu_d)$, define the determinant
\[
V(\mu_1,\ldots,\mu_d)=\left\vert\begin{matrix}
t_1^{\mu_1}&t_2^{\mu_1}&\cdots&t_d^{\mu_1}\\
t_1^{\mu_2}&t_2^{\mu_2}&\cdots&t_d^{\mu_2}\\
\vdots&\vdots&\ddots&\vdots\\
t_1^{\mu_d}&t_2^{\mu_d}&\cdots&t_d^{\mu_d}\\
\end{matrix}\right\vert.
\]
Then the Schur function is
\[
S_{\lambda}(t_1,\ldots,t_d)=
\frac{V(\lambda_1+d-1,\lambda_2+d-2,\ldots,\lambda_{d-1}+1,\lambda_d)}
{V(d-1,d-2,\ldots,1,0)}.
\]
The Schur functions play the role of
characters of the corresponding representation of $GL_d$.
{\it If we know the Hilbert series $H(C_{nd},t_1,\ldots,t_d)$,
then we can uniquely determine the $GL_d$-module structure of $C_{nd}$}.

\section{Concrete computations}

We give a survey of some concrete results about generators, defining
relations and Hilbert series of the algebras of matrix invariants.
We start with $2\times 2$ matrices.

The first important reduction is the following. We take the generic
$2\times 2$ matrix $X$ and present it in the form
\[
X=\left(\begin{matrix}
x_{11}&x_{12}\\
x_{21}&x_{22}\\
\end{matrix}\right)=
\left(\begin{matrix}
\frac{x_{11}+x_{22}}{2}&0\\
0&\frac{x_{11}+x_{22}}{2}\\
\end{matrix}\right)
+\left(\begin{matrix}
\frac{x_{11}-x_{22}}{2}&x_{12}\\
x_{21}&-\frac{x_{11}-x_{22}}{2}\\
\end{matrix}\right)
\]
\[
=\frac{\text{\rm tr}(X)}{2}E
+\left(\begin{matrix}
y_{11}&y_{12}\\
y_{21}&-y_{11}\\
\end{matrix}\right)
=\frac{\text{\rm tr}(X)}{2}E+Y,
\]
where $E$ is the identity matrix and
\[
Y=\left(\begin{matrix}
y_{11}&y_{12}\\
y_{21}&-y_{11}\\
\end{matrix}\right)
\]
is a generic $2\times 2$ {\it traceless} matrix.
This reduction implies that $C_{2d}$ is generated by
$\text{\rm tr}(X_i)$, $i=1,\ldots,d$, and
$\text{\rm tr}(Y_{i_1}\cdots Y_{i_k})$, where
$Y_1,\ldots,Y_d$ are generic traceless matrices
and $k\geq 2$. Since the class of nilpotency in the Nagata-Higman
theorem is $N(2)=3$ for $n=2$, we obtain that it is sufficient to consider the
cases $k=2,3$ only. Similarly, $T_{2d}$ is generated by
$C_{2d}$ and $Y_1,\ldots,Y_d$. It is also easy to see:

\begin{proposition}
The algebra $C_{2d}$ has the presentation
\[
C_{2d}=\left(K[\text{\rm tr}(X_1),\ldots,\text{\rm tr}(X_d)]\right)
\otimes_K C_0,
\]
where
\[
C_0=K[\text{\rm tr}(Y_iY_j),\text{\rm tr}(Y_pY_qY_r)\mid
1\leq i\leq j\leq d,1\leq p<q<r\leq d]/I
\]
is the algebra generated by products of generic traceless $2\times 2$ matrices.
Hence we may choose all defining relations as linear combinations of the traces
of products of traceless matrices.
\end{proposition}

The description of $T_{2d}$ is easier than that of $C_{2d}$.

\begin{theorem}
{\rm (i) (Procesi \cite{P3})} The noncommutative trace algebra $T_{2d}$ is
isomorphic to the tensor product of $K$-algebras
\[
K[\text{\rm tr}(X_1),\ldots,\text{\rm tr}(X_m)] \otimes_K W_d,
\]
where $W_d$ is the associative algebra generated by the
generic traceless matrices $Y_1,\ldots,Y_m$.

{\rm (ii) (Razmyslov \cite{R3})} The algebra of the generic traceless
matrices $W_d$ is isomorphic to the factor algebra
$K\langle y_1,\ldots,y_d\rangle/I(M_2(K),sl_2(K))$
of the free associative algebra $K\langle y_1,\ldots,y_d\rangle$,
where $sl_2(K)$ is the Lie algebra of all traceless
$2 \times 2$ matrices and
$I(M_2(K),sl_2(K))$ is the ideal of all polynomials
in $K\langle y_1,\ldots,y_d\rangle$ which vanish
under the substitutions of $y_i$ by elements in $sl_2(K)$.
{\rm (Such polynomials are called weak polynomial identities for the pair
$(M_2(K),sl_2(K))$.)}
The ideal $I(M_2(K),sl_2(K))$ is generated
as a weak T-ideal by the weak polynomial identity
$[x_1^2,x_2] = 0$, i.e. $I(M_2(K),sl_2(K))$
is the minimal ideal of weak
polynomial identities containing the element $[x_1^2,x_2]$.
An equivalent description is that as an ideal of the free algebra,
$I(M_2(K),sl_2(K))$ is generated by all elements
$[uv+vu,w]$, where $u,v,w$ are all possible commutators
$[y_{j_1},\ldots,y_{j_k}]$ in the variables $y_1,\ldots,y_d$.

{\rm (iii) (Drensky and Koshlukov \cite{DK})} The ideal
$I(M_2(K),sl_2(K))$ is the minimal ideal of
the free associative algebra
$K\langle y_1,\ldots,y_d\rangle$ which is invariant under the diagonal
action of $GL_d$ and contains the elements
$[y_1^2,y_2]$ and
\[
s_4(y_1,y_2,y_3,y_4)=\sum_{\sigma\in S_4}\text{\rm sign}(\sigma)
y_{\sigma(1)}y_{\sigma(2)}y_{\sigma(3)}y_{\sigma(4)},
\]
the second polynomial appears for $d \geq 4$
only. Hence the algebra of $2 \times 2$ generic traceless matrices
has a uniform set of defining relations for any $d \geq 2$.

{\rm (iv) (Procesi \cite{P3})} As a $GL_d$-module $W_d$ has the description
\[
W_d\cong \bigoplus W(\lambda_1,\lambda_2,\lambda_3),
\]
where the sum is on all partitions $\lambda$ in at most three parts.

{\rm (v) (Formanek \cite{F7})} The $GL_d$-module $T_{2d}$ has the decomposition
\[
T_{2d}\cong\bigoplus (\lambda_1-\lambda_2+1)(\lambda_2-\lambda_3+1)(\lambda_3-\lambda_4+1)
W(\lambda_1,\lambda_2,\lambda_3,\lambda_4).
\]

{\rm (vi) (See \cite{D2, D, F7, P3, LB1, LB2} for other descriptions of the Hilbert series.)}
The Hilbert series of $T_{2d}$ is
\[
H(T_{2d},t_1,\ldots,t_d)=\prod_{i=1}^d\frac{1}{1-t_i}\sum S_{(\lambda_1,\lambda_2,\lambda_3)}(t_1,\ldots,t_d).
\]
\end{theorem}

Now we give some results on $C_{2d}$. We follow the way we used in the previous theorem.
For more details, especially for the Hilbert series and the $GL_d$-module decomposition
of $C_{2d}$, see \cite{F7, LB2, P3} or \cite{DF}.

\begin{theorem}
{\rm (i)} The commutative trace algebra $C_{2d}$ is
isomorphic to the tensor product of $K$-algebras
\[
K[\text{\rm tr}(X_1),\ldots,\text{\rm tr}(X_m)] \otimes_K C(W_d),
\]
where $C(W_d)$ is the centre of the algebra $W_d$ defined above.

{\rm (ii)} As a subalgebra of $W_d$, its centre $C(W_d)$ is generated by
\[
Y_i^2,\quad i=1,\ldots,d, \quad Y_iY_j+Y_jY_i,\quad 1\leq i<j\leq d,
\]
\[
s_3(Y_i,Y_j,Y_k)=\sum_{\sigma\in S_3}\text{\rm sign}(\sigma)
Y_{\sigma(1)}Y_{\sigma(2)}Y_{\sigma(3)},\quad 1\leq i<j<k\leq d,
\]
where the symmetric group $S_3$ acts on $\{i,j,k\}$.

{\rm (iii)} As a $GL_d$-module $W_d$ has the decomposition
\[
W_d\cong \bigoplus W(\lambda_1,\lambda_2,\lambda_3),
\]
where the sum is on all partitions $\lambda$ in at most three parts such that
both $\lambda_1-\lambda_2$ and $\lambda_2-\lambda_3$ are even.
\end{theorem}

Concerning the defining relations of $C_{2d}$, the case $d=2$ is trivial.
Formanek, Halpin and Li \cite{FHL} showed that
$C_{22}$ is generated by the algebraically independent elements
\[
\text{\rm tr}(X_1), \text{\rm tr}(X_2), \text{\rm det}(X_1), \text{\rm det}(X_2), \text{\rm tr}(X_1X_2).
\]
For $d=3$ Sibirskii \cite{Si} found one relation between the generators of $C_{23}$
and, using the Hilbert series of $C_{23}$, Formanek \cite{F7} proved that there are no more relations.
In the general case, the description of the defining relations of $C_{2d}$ is reduced to a similar description
of the defining relations of the subalgebra of $C_{2d}$ generated by
\[
\text{\rm tr}(Y_i^2),\quad i=1,\ldots,d,\quad \text{\rm tr}(Y_iY_j),\quad 1\leq i<j\leq d,
\]
\[
\text{\rm tr}(Y_iY_jY_k),\quad 1\leq i<j<k\leq d,
\]
where $Y_1,\ldots,Y_d$ are generic traceless $2\times 2$ matrices. Since $GL_2$ acts on
the generic matrices by conjugation, we may replace its action with the action of $SL_2$ and even with
the action of $PSL_2$. Since $PSL_2({\mathbb C})\cong SO_3({\mathbb C})$,
the special orthogonal group, we may apply invariant theory of
special linear groups. (The restriction $K=\mathbb C$ is not essential in the final version of the result.)

We consider the action of the special orthogonal group $SO_3=SO_3(K)$, i.e. the group of orthogonal $3\times 3$ matrices
with determinant 1, on the polynomial algebra in $3d$ variables
\[
K[u_i^{(1)},u_i^{(2)},u_i^{(3)}\mid i=1,\ldots, d],
\]
induced by the action of $SO_3$ on the three-dimensional
vectors $u_i=(u_i^{(1)},u_i^{(2)},u_i^{(3)})$. It is a classical result
that the algebra of invariants $K[u_i^{(j)}]^{SO_3}$
is generated by all scalar products
\[
\langle u_i,u_j\rangle=u_i^{(1)}u_j^{(1)}+u_i^{(2)}u_j^{(2)}+u_i^{(3)}u_j^{(3)},
\quad 1\leq i\leq j\leq d,
\]
and all $3\times 3$ determinants of the coordinates
\[
\Delta(u_i,u_j,u_k)=\left\vert\begin{matrix}
u_i^{(1)}&u_j^{(1)}&u_k^{(1)}\\
u_i^{(2)}&u_j^{(2)}&u_k^{(2)}\\
u_i^{(3)}&u_j^{(3)}&u_k^{(3)}\\
\end{matrix}\right\vert, \quad 1\leq i<j<k\leq d.
\]
The defining relations express the fact that the underlying vector space is three-dimensional and every
four vectors are linearly dependent. In particular, they use the properties of the Gram determinant:
\[
\Gamma_4(u_i,u_j,u_k,u_l;u_p,u_q,u_r,u_s)=\left\vert\begin{matrix}
\langle u_i,u_p\rangle&\langle u_i,u_q\rangle&\langle u_i,u_r\rangle&\langle u_i,u_s\rangle\\
\langle u_j,u_p\rangle&\langle u_j,u_q\rangle&\langle u_j,u_r\rangle&\langle u_j,u_s\rangle\\
\langle u_k,u_p\rangle&\langle u_k,u_q\rangle&\langle u_k,u_r\rangle&\langle u_k,u_s\rangle\\
\langle u_l,u_p\rangle&\langle u_l,u_q\rangle&\langle u_l,u_r\rangle&\langle u_l,u_s\rangle\\
\end{matrix}\right\vert=0,
\]
\[
1\leq i<j<k<l\leq d,\quad 1\leq p<q<r<s\leq d,
\]
\[
\Delta(u_i,u_j,u_k)\Delta(u_p,u_q,u_r)-\Gamma_3(u_i,u_j,u_k;u_p,u_q,u_r)=0,
\]
\[
\langle u_p,u_i\rangle\Delta(u_j,u_k,u_l)
-\langle u_p,u_j\rangle\Delta(u_i,u_k,u_l)
\]
\[
+\langle u_p,u_k\rangle\Delta(u_i,u_j,u_l)
-\langle u_p,u_l\rangle\Delta(u_i,u_j,u_k)=0.
\]
In order to apply invariant theory of $SO_3$ we need a scalar product (i.e. nondegenerate
symmetric bilinear form) on $sl_2(K)$. We use the trace and define
\[
\langle u,v\rangle=\text{\rm tr}(uv),\quad u,v\in sl_2(K).
\]
The following result gives the generators and the defining relations of $C_{2d}$ for $d\geq 2$.
It is a translation of the description of the invariants of $SO_3$.

\begin{theorem}
{\rm (i)} The algebra $C_{2d}$ is generated by
\[
\text{\rm tr}(X_i),\quad \text{\rm tr}(Y_i^2),\quad \text{\rm tr}(Y_iY_j),\quad
\text{\rm tr}(s_3(Y_i,Y_j,Y_k)),
\]
where $i,j,k=1,\ldots,d$, and in the traces involving two or three traceless matrices we require
$i<j$ or $i<j<k$, respectively.

{\rm (ii) Drensky \cite{D1})} The defining relations of $C_{2d}$ with respect to the above generators are
\[
\text{\rm tr}(s_3(Y_i,Y_j,Y_k))\text{\rm tr}(s_3(Y_p,Y_q,Y_r))
+18\left\vert\begin{matrix}
\text{\rm tr}(Y_iY_p)&\text{\rm tr}(Y_iY_q)&\text{\rm tr}(Y_iY_r)\\
\text{\rm tr}(Y_jY_p)&\text{\rm tr}(Y_jY_q)&\text{\rm tr}(Y_jY_r)\\
\text{\rm tr}(Y_kY_p)&\text{\rm tr}(Y_kY_q)&\text{\rm tr}(Y_kY_r)\\
\end{matrix}\right\vert=0,
\]
\[
\text{\rm tr}(Y_pY_i)\text{\rm tr}(s_3(Y_j,Y_k,Y_l))
-\text{\rm tr}(Y_pY_j)\text{\rm tr}(s_3(Y_i,Y_k,Y_l))
\]
\[
+\text{\rm tr}(Y_pY_k)\text{\rm tr}(s_3(Y_i,Y_j,Y_l))
-\text{\rm tr}(Y_pY_l)\text{\rm tr}(s_3(Y_i,Y_j,Y_k))=0,
\]
where, again, $i,j,k,p,q,r=1,\ldots,d$, and, where necessary, we require $i<j<k<l$ and $p<q<r$.
\end{theorem}

In order to work efficiently with an algebra $R=K[x_1,\ldots,x_p]/I$,
it is not sufficient to know the generators of the ideal $I$. For computational purposes one needs also
the {\it Gr\"obner basis} of $I$ with respect to some ordering on the monomials of $K[x_1,\ldots,x_p]$, see e.g.
the book by Adams and Loustaunau \cite{AL}. The Gr\"obner basis of $C_{2d}$ is given by Domokos and Drensky \cite{DD},
see their paper for more details.

Now we consider two generic $3\times 3$ matrices.
Using the Molien-Weyl formula, Teranishi \cite{T1} calculated the Hilbert series of $C_{32}$, namely,
\[
H(C_{32},t_1,t_2)=
\frac{1+t_1^3t_2^3}
{(1-t_1)(1-t_2)q_2(t_1,t_2)q_3(t_1,t_2)(1-t_1^2t_2^2)},
\]
where
\[
q_2(t_1,t_2)=(1-t_1^2)(1-t_1t_2)(1-t_2^2),
\]
\[
q_3(t_1,t_2)=(1-t_1^3)(1-t_1^2t_2)(1-t_1t_2^2)(1-t_2^3).
\]
He also found the following system of generators of $C_{32}$:
\[
\begin{array}{c}
\text{\rm tr}(X_1),\text{\rm tr}(X_2),\text{\rm tr}(X_1^2),\text{\rm tr}(X_1X_2),\text{\rm tr}(X_2^2), \\
\\
\text{\rm tr}(X_1^3),\text{\rm tr}(X_1^2X_2),\text{\rm tr}(X_1X_2^2),\text{\rm tr}(X_2^3),
\text{\rm tr}(X_1^2X_2^2),\text{\rm tr}(X_1^2X_2^2X_1X_2),
\end{array}
\]
where $X_1,X_2$ are generic $3 \times 3$ matrices.
He showed that the first ten of these generators
form a homogeneous system of parameters of $C_{32}$ and
$C_{32}$ is a free module with generators 1 and $\text{\rm tr}(X_1^2X_2^2X_1X_2)$
over the polynomial algebra on these ten elements.

Abeasis and Pittaluga \cite{AP} found a system of generators of
$C_{3d}$ in terms of representation theory of the symmetric and general linear groups,
in the spirit of its use in theory of PI-algebras.
They showed that $C_{3d}$ has a minimal system of generators which spans
a $GL_d$-module isomorphic to
\[
G=W(1)\oplus W(2)\oplus W(3)\oplus W(1^3)\oplus W(2^2)\oplus W(2,1^2)
\]
\[
\oplus W(3,1^2)\oplus W(2^2,1)\oplus W(1^5)\oplus W(3^2)\oplus W(3,1^3).
\]
(The partitions in \cite{AP} are given in ``Francophone'' way, i.e., transposed to ours.)

It follows from the description of the generators of $C_{32}$ given by Teranishi \cite{T1}, that
$\text{\rm tr}(X_1^2X_2^2X_1X_2)$ satisfies a quadratic equation with coefficients
depending on the other ten generators. The explicit (but very complicated)
form of the equation was found by Nakamoto \cite{N}, over $\mathbb Z$,
with respect to a slightly different system of generators.
A much simpler description of $C_{32}$ was obtained by Aslaksen, Drensky, and Sadikova \cite{ADS}.
The following generators are in the spirit of the ideas of \cite{AP}.

\begin{proposition}
Let $X_1,X_2$ and $Y_1,Y_2$ be, respectively, two generic and two generic traceless
$3\times 3$ matrices.
The algebra $C_{32}$ is generated by
\[
\begin{array}{c}
\text{\rm tr}(X_1),\text{\rm tr}(X_2),\text{\rm tr}(Y_1^2),\text{\rm tr}(Y_1Y_2),\text{\rm tr}(Y_2^2), \\
\\
\text{\rm tr}(Y_1^3),\text{\rm tr}(Y_1^2Y_2),\text{\rm tr}(Y_1Y_2^2),\text{\rm tr}(Y_2^3),\\
\\
V=\text{\rm tr}(Y_1^2Y_2^2)-\text{\rm tr}(Y_1Y_2Y_1Y_2),\quad
W=\text{\rm tr}(Y_1^2Y_2^2Y_1Y_2)-\text{\rm tr}(Y_2^2Y_1^2Y_2Y_1).
\end{array}
\]
\end{proposition}

Now we define the following elements of $C_{32}$:
\[
U=\left| \begin{array}{cc}
\text{\rm tr}(Y_1^2)& \text{\rm tr}(Y_1Y_2) \\
\text{\rm tr}(Y_1Y_2) & \text{\rm tr}(Y_2^2) \\
\end{array} \right|,
\]
\[
W_1=U^3,  \quad
W_2=U^2V, \quad
W_4=UV^2, \quad
W_7=V^3,
\]
\[
W_5=V\left| \begin{array}{ccc}
\text{\rm tr}(Y_1^2)& \text{\rm tr}(Y_1Y_2) & \text{\rm tr}(Y_2^2) \\
\text{\rm tr}(Y_1^3)& \text{\rm tr}(Y_1^2Y_2) & \text{\rm tr}(Y_1Y_2^2) \\
\text{\rm tr}(Y_1^2Y_2) & \text{\rm tr}(Y_1Y_2^2) & \text{\rm tr}(Y_2^3) \\
\end{array} \right| ,
\]
\[
W_6=\left| \begin{array}{cc}
\text{\rm tr}(Y_1^3)& \text{\rm tr}(Y_1Y_2^2) \\
\text{\rm tr}(Y_1^2Y_2) & \text{\rm tr}(Y_2^3) \\
\end{array} \right|^2
-4\left| \begin{array}{cc}
\text{\rm tr}(Y_2^3)& \text{\rm tr}(Y_1Y_2^2) \\
\text{\rm tr}(Y_1Y_2^2) & \text{\rm tr}(Y_1^2Y_2) \\
\end{array} \right|\left| \begin{array}{cc}
\text{\rm tr}(Y_1^3)& \text{\rm tr}(Y_1^2Y_2) \\
\text{\rm tr}(Y_1^2Y_2) & \text{\rm tr}(Y_1Y_2^2) \\
\end{array} \right|,
\]
\[
W_3'=U\left| \begin{array}{ccc}
\text{\rm tr}(Y_1^2)& \text{\rm tr}(Y_1Y_2) & \text{\rm tr}(Y_2^2) \\
\text{\rm tr}(Y_1^3)& \text{\rm tr}(Y_1^2Y_2) & \text{\rm tr}(Y_1Y_2^2) \\
\text{\rm tr}(Y_1^2Y_2) & \text{\rm tr}(Y_1Y_2^2) & \text{\rm tr}(Y_2^3) \\
\end{array} \right| ,
\]
where $U,V$ are defined above. Finally, we define one more element $W_3''$ as follows.
Recall that a linear mapping $\delta$ of an algebra $R$ is a {\it derivation}
if $\delta(rs)=\delta(r)s+r\delta(s)$ for all $r,s\in R$. We consider the derivation
$\delta$ of $C_{32}$ which commutes with the trace and satisfies the conditions
\[
\delta(X_1)=0,\quad \delta(X_2)=X_1,\quad \delta(Y_1)=0,\quad \delta(Y_2)=Y_1.
\]
Then
\[
W_3''=\frac{1}{144}\sum_{i=0}^6(-1)^i\delta^i(\text{\rm tr}^3(Y_2^2))
\delta^{6-i}(\text{\rm tr}^2(Y_2^3)).
\]

The following theorem gives the defining relation of $C_{32}$.
It uses representation theory of $GL_2$, combinatorics, computations by hand and
easy computer calculations with standard functions of
Maple.

\begin{theorem}
{\rm (Aslaksen, Drensky, Sadikova \cite{ADS})}
The algebra of invariants $C_{32}$ of two $3\times 3$ matrices
is generated by the elements from the previous theorem,
subject to the defining relation
\[
W^2-\left( \frac{1}{27}W_1-\frac{2}{9}W_2+\frac{4}{15}W_3'+\frac{1}{90}W_3''
+\frac{1}{3}W_4-\frac{2}{3}W_5-\frac{1}{3}W_6-\frac{4}{27}W_7 \right)=0.
\]
\end{theorem}

The calculation of the Hilbert series of $C_{nd}$ and $T_{nd}$
based on the Molien-Weyl formula is quite complicated because requires
evaluations of multiple integrals. Van den Bergh \cite{VB2}
sujected a way which involves graph theory. As a consequence, he established
important properties of $H(C_{2d},t_1,\ldots,t_d)$ and $H(T_{2d},t_1,\ldots,t_d)$.

Berele and Stembridge \cite{BS} applied the method of van den Bergh \cite{VB2} and calculated the Hilbert series of $T_{32}$.
Using the above results of Aslaksen, Drensky, and Sadikova on $C_{32}$
and the explicit form of the Hilbert series of $T_{32}$, Benanti and Drensky \cite{BD}
found a polynomial subalgebra $S$ of $C_{32}$ and a finite set of generators of the free $S$-module
$T_{32}$. They gave also a set of defining relations of $T_{32}$ as an algebra
and a Gr\"obner basis of the corresponding ideal. (See the survey article by Ufnarovski \cite{U}
for a background on Gr\"obner bases in
the noncommutative case, as well as the paper by Mikhalev and Zolotykh \cite{MZ} which is closer to
the situation in \cite{BD}.)

For two generic $4\times 4$ matrices, the Hilbert series of $C_{42}$ was calculated (with some typos)
by Teranishi \cite{T1, T2} and corrected by Berele and Stembridge \cite{BS}.
Teranishi found also a homogeneous system of parameters and a system of generators,
in the spirit of the $3\times 3$ case.
Recently, Drensky and Sadikova \cite{DS} have found another system of generators of $C_{42}$
which is minimal and seems to be more convenient for concrete calculations.

The Hilbert series of a graded vector space with $GL_d$-module structure determines uniquely
its decomposition into irreducible submodules. Hence, in principle, one may calculate the multiplicities
$m(\lambda)$ if one knows the concrete form of the Hilbert series. Berele \cite{B}
used the Hilbert series of $C_{32}$ found by Teranishi \cite{T1} and described the asymptotics of
$m(\lambda_1,\lambda_2)$. (Due to a technical error (an omitted summand) some of the coefficients
of the polynomials in the asymptotics of Berele
are slightly different from the real ones.) Another approach to the problem was suggested by Drensky and Genov
\cite{DG1}. Let
\[
f(t_1,t_2)=\sum_{i,j\geq 0}a_{ij}t_1^it_2^j,
\]
$a_{ij}\in K$, $a_{ij}=a_{ji}$,
be a symmetric function in two variables which is a formal power series from $K[[t_1,t_2]]$.
We present it in the form
\[
f(t_1,t_2)=\sum_{\lambda_1\geq\lambda_2}m(\lambda_1,\lambda_2)S_{(\lambda_1,\lambda_2)}(t_1,t_2)
\]
and want to find the multiplicities $m(\lambda_1,\lambda_2)$.
In most of the cases which we consider, $f(t_1,t_2)$ is given explicitly as a rational function.
So, it is natural to express $m(\lambda_1,\lambda_2)$ not in terms of the coefficients $a_{ij}$
but in a more direct way. We introduce the generating function of the multiplicities
\[
M(f,t,u)=\sum_{\lambda_1\geq\lambda_2}m(\lambda_1,\lambda_2)t^{\lambda_1}u^{\lambda_2}\in K[[t,u]]
\]
and call it the {\it multiplicity series} of $f(t_1,t_2)$. It is more convenient to introduce
a new variable $v=tu$ and to consider the series
\[
M'(f,t,v)=\sum_{\lambda_1\geq\lambda_2}m(\lambda_1,\lambda_2)t^{\lambda_1-\lambda_2}v^{\lambda_2}\in K[[t,v]],
\]
because the mapping $M': K[[t_1,t_2]]^{S_2}\to K[[t,v]]$ is a bijective linear mapping which
is continuous with respect to the {\it formal power series topology}. It is easy to see that
$f(t_1,t_2)$ and $M'(f,t,v)$ are related by
\[
f(t_1,t_2)=\frac{t_1M'(f,t_1,t_1t_2)-t_2M'(f,t_2,t_1t_2)}{t_1-t_2}.
\]
Hence, if we have a potential candidate $h(t,v)$ for $M'(f,t,v)$, it is easy
to verify whether $h(t,v)=M'(f,t,v)$.

Also, the elementary symmetric function $e_2=t_1t_2$ behaves like a constant,
\[
M'(g(t_1t_2)f(t_1,t_2),t,v)=g(v)M'(f,t,v),
\]
and this simplifies the calculations. Applying quite complicated (also technically) arguments,
Drensky and Genov \cite{DG1} found the multiplicity series of the Hilbert series of $C_{32}$.
They also corrected the technical errors in \cite{B}.

\begin{theorem} {\rm (i)} \cite{DG1} The multiplicity series of the
Hilbert series of the algebra $C_{32}$ of invariants of two
$3\times 3$ matrices is
\[
M'(H(C_{32},t_1,t_2),t,v)=
\frac{1}{(1-v^2)(1-v^3)^2}\times
\]
\[
\times\left(
\frac{(1+v^2+v^4)((1+v^2)(1-t^2v)+2tv(1-v))}
{3(1-v)(1-v^2)^3(1-t)^2(1-t^2)}+\right.
\]
\[
+\frac{(1-v)(1+tv)}{3(1-v^2)(1-t)(1-t^2)}
+\frac{(1-v^2)(1-tv)}{3(1-v^3)(1-t^3)}
\]
\[
\left. -\frac{v^3((1-v+v^2)(1-t^2v^2)+tv(1-v^2))}
{(1-v)(1-v^2)^2(1-v^4)(1-t)(1-t^2)(1-tv)}\right).
\]

{\rm (ii)} \cite{B, DG1} Let $\lambda=(p,q)$ and let
$m(p,q)$ be the multiplicity of $S_{(p,q)}(t_1,t_2)$ in $H(C_{32},t_1,t_2)$.
Then for $p>2q\geq 0$
\[
m(p,q)=\frac{q^7}{7!2^5.3^2}+\frac{(p-q)q^6}{6!2^4.3^2}
+\frac{(p-q)^2q^5}{2!5!2^33^2}
+{\mathcal O}((p+q)^6)
\]
\[
=\frac{p^2q^5}{17280}-\frac{11pq^6}{103680}+\frac{71q^7}{1451520}+
{\mathcal O}((p+q)^6);
$$
for $2q\geq p\geq q\geq 0$
$$
m(p,q)=\frac{q^7}{7!2^5.3^2}+\frac{(p-q)q^6}{6!2^4.3^2}
+\frac{(p-q)^2q^5}{2!5!2^33^2}-\frac{(2q-p)^7}{7!2^5.3^2}
+{\mathcal O}((p+q)^6)
\]
\[
=\frac{p^7}{1451520}-\frac{p^6q}{103680}+\frac{p^5q^2}{17280}-\frac{p^4q^3}{5184}
+\frac{p^3q^4}{2592}-\frac{7p^2q^5}{17280}+
\frac{7pq^6}{34560}-\frac{19q^7}{483840}+{\mathcal O}((p+q)^6).
\]
\end{theorem}

Later the methods for calculating the multiplicity series of symmetric functions of special
kinds were significantly improved \cite{DG2}. The Hilbert series of $C_{42}$
calculated by Teranishi \cite{T1, T2} (with some typos corrected in \cite{BS})
and the Hilbert series of $T_{32}$ and $T_{42}$ calculated by Berele and Stembridge \cite{BS}
allowed to express their multiplicity series and to determine the asympotics of the multiplicities.
We shall state simplified versions of the results:
\begin{theorem} {\rm (i) (Drensky, Genov, Valenti \cite{DGV})}
The multiplicities $m_{(\lambda_1,\lambda_2)}(C_{32})$ and
$m_{(\lambda_1,\lambda_2)}(T_{32})$ of the Hilbert series of $C_{32}$ and $T_{32}$, respectively, are related by
\[
m_{(\lambda_1,\lambda_2)}(T_{32})\approx
9m_{(\lambda_1,\lambda_2)}(C_{32}).
\]

{\rm (ii)( Drensky and Genov \cite{DGV})}
Let $\lambda=(\lambda_1,\lambda_2)$.
The multiplicities $m_{\lambda}(C_{42})$ of the Hilbert series of
$C_{42}$ satisfy the condition
\[
m_{\lambda}(C_{42})=\begin{cases}
m_1+{\mathcal O}((\lambda_1+\lambda_2)^{13}), & \text{if $\lambda_1>3\lambda_2$,}\\
m_1+m_2+{\mathcal O}((\lambda_1+\lambda_2)^{13}), & \text{if $3\lambda_2\geq\lambda_1>2\lambda_2$,}\\
m_1+m_2+m_3+{\mathcal O}((\lambda_1+\lambda_2)^{13}), & \text{if $2\lambda_2\geq \lambda_1$,}\\
\end{cases}
\]
where
\[
m_1=\frac{(\lambda_1-\lambda_2)^3\lambda_2^{11}}{11!3!2^83^2}
-\frac{(\lambda_1-\lambda_2)^2\lambda_2^{12}}{12!2!2^83^3}
+\frac{127(\lambda_1-\lambda_2)\lambda_2^{13}}{13!2^{10}3^4}
-\frac{305\lambda_2^{14}}{14!2^93^5},
\]
\[
m_2=
\frac{(3\lambda_2-\lambda_1)^{14}}{14!2^{10}3^55^2},
\]
\[
m_3=-\frac{(\lambda_1-\lambda_2)(2\lambda_2-\lambda_1)^{13}}{13!2^{10}3^25}
-\frac{7(2\lambda_2-\lambda_1)^{14}}{14!2^93\cdot 5^2}.
\]
The multiplicities $m_{\lambda}(T_{42})$ satisfy
\[
m_{\lambda}(T_{42})=16m_{\lambda}(C_{42})+{\mathcal O}((\lambda_1+\lambda_2)^{13}).
\]
\end{theorem}

We want to mention that Berele and Stembridge \cite{BS} computed also the Hilbert series of
$C_{33}$ and $T_{33}$ but the methods of \cite{DG1, DG2, DG3, DGV} do not work successfully
for symmetric functions in three variables. One can introduce the multiplicity series
of a symmetric function in any number of variables, generalizing in an obvious way the case
of symmetric functions in two variables. A recent theorem of Berele \cite{B2} gives
the rationality of the multiplicity series
of a class of rational symmetric functions
in any number of variables, including the Hilbert series of $C_{nd}$ and $T_{nd}$.
Unfortunately, it is not clear how to perform the concrete calculations, even for the Hilbert series
of $C_{33}$ and $T_{33}$.

\end{document}